\documentclass{article}
\usepackage[a4paper,mag=1000,includefoot,left=2cm,right=2cm,top=2cm,bottom=2cm,headsep=1cm,footskip=1cm]{geometry}
\usepackage{pgfplotstable,booktabs,colortbl}
\usepackage{graphics}
\usepackage{epsfig}
\usepackage{amssymb}

\usepackage[english,english]{babel}
\usepackage[T2A]{fontenc}
\usepackage[cp1251]{inputenc}
\usepackage{amssymb,amsmath,indentfirst,bm}
\newtheorem{theorem}{Theorem}[section]
\newtheorem{lemma}[theorem]{Lemma}

 \renewcommand{\geq}{\geqslant}
\begin{document}
\large
\centerline{\Large A DIRECT FAST FFT-BASED IMPLEMENTATION}
\smallskip
\centerline{\Large FOR HIGH ORDER FINITE ELEMENT METHOD}
\smallskip
\centerline{\Large ON RECTANGULAR PARALLELEPIPEDS FOR PDE}
% A direct fast FFT-based implementation for high order finite element method on rectangular parallelepipeds for PDE
\medskip\par\centerline{\large A. ZLOTNIK
\footnote{National Research University Higher School of Economics, Myasnitskaya 20, 101000 Moscow, Russia
(azlotnik2008@gmail.com)}
and I. ZLOTNIK
\footnote{Settlement Depository Company, 2-oi Verkhnii Mikhailovskii proezd 9, building 2, 115419 Moscow, Russia
(ilya.zlotnik@gmail.com)}}
\medskip
\begin{abstract}
We present a new direct logarithmically optimal in theory and fast in practice algorithm to implement the high order finite element method on multi-dimensional rectangular parallelepipeds for solving PDEs of the Poisson kind.
The key points are the fast direct and inverse FFT-based algorithms for decomposition in eigenvectors of the 1D eigenvalue problems for the high order FEM.
The algorithm can further be used for numerous applications, in particular, to implement the high order finite element methods for various time-dependent PDEs.
\end{abstract}
\smallskip\textbf{Key words.} Fast direct algorithm, high order finite element method, FFT, Poisson equation.
\medskip\par
\textbf{AMS subject classifications.} 65F05, 65F15, 65M60, 65T99.

\section{Introduction}
\label{}

We present new direct fast algorithm to implement $n$th order ($n\geq 2$) finite element method (FEM) on rectangular parallelepipeds \cite{C02} for solving $N$-dimensional PDEs, $N\geq 2$, like the Poisson one with the Dirichlet boundary condition.
The algorithm generalizes the well-known one in the case of the bilinear elements ($n=1$) or standard finite-difference schemes \cite{BFK11,SN78,S77} and utilizes the discrete fast Fourier transforms (FFTs) \cite{BRY07}.
The key points are the fast direct and inverse algorithms for decomposition in eigenvectors of the 1D eigenvalue problems for the high order FEM; this solves the known problem, see \cite[p. 271]{BFK11}.
The algorithm is logarithmically optimal with respect to the number of elements.
It also demonstrates rather mild growth in $n$ starting from the known case $n=1$ and is fast in practice, for example, the 2D FEM system for $2^{20}$ elements of the 9th order containing almost $85\cdot10^6$ unknowns is solved in less than 2 min on an ordinary laptop, see Fig. \ref{fig:TIME:2D:ASUS} below.
The algorithm can further serve for a variety of applications including general 2nd order elliptic equations (as a preconditioner), for  the $N$-dimensional heat, wave or time-dependent Schr\"{o}dinger PDEs.
It can be applied for some non-rectangular domains, in particular, by involving meshes topologically equivalent to rectangular ones \cite{D96}.
Other standard boundary conditions can be covered as well \cite{ZZDAN16}; moreover,
the structure of the algorithm is valuable for wave problems with non-local boundary conditions, see \cite{BFK11,DZ06,DZ07,ZZ12}, whence our own interest arose.
The algorithm is also highly parallelizable.

\section{Algorithms}
\label{A}

\smallskip\par 1. We first need to consider in detail the FEM for the simplest 1D eigenvalue ODE problem
\begin{gather}
 -u''(x)=\lambda u(x)\ \ \text{on}\ \ [0,X],\ \ u(0)=u(X)=0,\ \ u(x)\not\equiv 0.
\label{eq:diff eig pr}
\end{gather}
\par We introduce the uniform mesh with the nodes $x_j=jh$, $j=\overline{0,K}$ (i.e., $0\leqslant j\leqslant K$)
and the step $h=X/K$.
Let $H_h^{(n)}[0,X]$ be the space of the piecewise-polynomial functions $\varphi\in C[0,X]$ such that
$\varphi(x)\in \mathcal{P}_n$ for $x\in [x_{j-1}, x_j]$, $j=\overline{1,K}$, with $\varphi(0)=\varphi(X)=0$;
here $\mathcal{P}_n$ is the space of polynomials having at most $n$th degree, $n\geq 2$.
\par Let $S_K^{(n)}$ be the space of vector functions $w$ such that
$w_j\in \mathbb{R}$ for $j=\overline{0,K}$ with $w_0=w_K=0$ and $w_{j-1/2}\in\mathbb{R}^{n-1}$, $j=\overline{1,K}$.
Clearly $\dim S_K^{(n)}=nK-1$.
A function $\varphi\in H_h^{(n)}[0,X]$ is uniquely defined by its values at the mesh nodes $\varphi_j=\varphi(x_j)$, $j=\overline{0,K}$, with $\varphi_0=\varphi_K=0$, and inside the elements $\varphi_{j-1/2}=\{\varphi(x_{j-1}+(l/n)h)\}_{l=1}^{n-1}$, $j=\overline{1,K}$, that form the element in $S_K^{(n)}$.
\par We utilize the following scaled operator form of the standard FEM discretization for problem \eqref{eq:diff eig pr}
\begin{gather}
 \mathcal{A}v=\lambda\mathcal{C}v,\ \ v\in S_K^{(n)},\ \ v\neq 0.
\label{eq:eig_glob}
\end{gather}
Here $\mathcal{A}=\mathcal{A}^T>0$ and $\mathcal{C}=\mathcal{C}^T>0$ are the global (scaled) stiffness and mass operators (matrices) acting in $S_K^{(n)}$ and together with $\lambda$ \textit{independent on} $h$; the true approximate eigenvalues are $\lambda_h=4h^{-2}\lambda$.
\par Let $A=\{A_{kl}\}_{k,l=0}^n$ and
$C=\{C_{kl}\}_{k,l=0}^n$ be the local stiffness and mass matrices related to the reference element $\sigma_0=[-1,1]$
with the following entries
\[
 A_{kl} = \int\nolimits_{\sigma_0}e'_k(x)e'_l(x)\,dx,\ \ C_{kl} = \int\nolimits_{\sigma_0}e_k(x)e_l(x)\,dx,
\]
where $\{e_l\}_{l=0}^n$ is the Lagrange basis in $\mathcal{P}_n$ such that
$e_l\bigl(-1+(2k)/n\bigr)=\delta_{kl}$, for $k,l=\overline{0,n}$,
and $\delta_{kl}$ is the Kronecker delta.
The matrices $A$, $C$ and the related matrix pencil have the following $3\times3$--block form
\begin{equation}
 A=
\left(\hspace{-4pt}
\begin{array}{ccc}
 a_0& a^T           & a_n\\
 a  & \widetilde{A} &\check{a}\\
 a_n& \check{a}^T   &a_0
\end{array}
\hspace{-4pt}\right),\ \
 C=
\left(\hspace{-4pt}
\begin{array}{ccc}
 c_0& c^T           & c_n\\
 c  & \widetilde{C} &\check{c}\\
 c_n& \check{c}^T   &c_0
\end{array}
\hspace{-4pt}\right),\ \
 G(\lambda):=A-\lambda C=\left(\hspace{-4pt}
\begin{array}{ccc}
 g_0(\lambda)& g^T(\lambda)          & g_n(\lambda)\\
 g(\lambda)  & \widetilde{G}(\lambda)&\check{g}(\lambda)\\
 g_n(\lambda)& \check{g}^T(\lambda)  &g_0(\lambda)
\end{array}
\hspace{-4pt}\right).
\label{eq:matrAC}
\end{equation}
Here $\widetilde{A}$, $\widetilde{C}$ and $\widetilde{G}(\lambda)=\widetilde{A}-\lambda \widetilde{C}$ are square matrices of order $n-1$ and $a,c,g(\lambda)=a-\lambda c\in\mathbb{R}^{n-1}$
whereas $\check{p}_l\equiv (Pp)_l=p_{n-l}$, $l=\overline{1,n-1}$, for $p\in\mathbb{R}^{n-1}$.
Let $\mathbb{R}_e^{n-1}$ and $\mathbb{R}_o^{n-1}$ be the subspaces of even and odd vectors in $\mathbb{R}^{n-1}$, i.e. such that $Pp=p$ and $Pp=-p$.
Clearly
$p=p_e+p_o$ with $p_e:=(p+\check{p})/2$ and $p_o:=(p-\check{p})/2$ and thus $\mathbb{R}^{n-1}=\mathbb{R}_e^{n-1}\oplus\mathbb{R}_o^{n-1}$ for $n\geq 3$; note that $\mathbb{R}_o^{n-1}=\{0\}$ for $n=2$.

Then problem \eqref{eq:eig_glob} can be represented in the following explicit form
\begin{gather*}
 g_n(\lambda) v_{j-1}+\check{g}(\lambda)\cdot v_{j-1/2}+2g_0(\lambda) v_j+g(\lambda)\cdot v_{j+1/2}+g_n(\lambda) v_{j+1}=0,\ \
 j=\overline{1,K-1},
\label{eq:eigpr1}
\\[1mm]
 g(\lambda)v_{j-1}+\widetilde{G}(\lambda)v_{j-1/2}+\check{g}(\lambda)v_j=0,\ \ j=\overline{1,K},
\label{eq:eigpr2}
\end{gather*}
with $v_0=v_K=0$, $v\not\equiv 0$;
see the similar problem for $\lambda\in \mathbb{C}$ on the uniform mesh on $[0,\infty)$ in \cite{ZZ12}.
Hereafter the symbol $\cdot$ denotes the inner product of vectors in $\mathbb{R}^{n-1}$.

\par We also consider the auxiliary eigenvalue problems on and inside the reference element $\sigma_0$
\begin{gather}
 Ae=\lambda Ce,\ \ e\in\mathbb{R}^{n+1},\ \ e\neq 0;\ \
 \widetilde{A}e=\lambda\widetilde{C}e,\ \ e\in\mathbb{R}^{n-1},\ \ e\neq 0,
\label{eq:eig_elem}
\end{gather}
where clearly $A\geq 0$, $C>0$ and $\widetilde{A}=\widetilde{A}^T>0$, $\widetilde{C}=\widetilde{C}^T>0$; see some their properties in  \cite{ZZ12}.
Denote by $S_n$ and $\tilde{S}_n$ their spectra.
Let $\{\lambda_0^{(l)},e^{(l)}\}_{l=1}^{n-1}$ be eigenpairs of the second problem \eqref{eq:eig_elem}.
\begin{lemma}
\label{prop1}
1. Any eigenvalue $\lambda_0^{(l)}$ is positive and at most double.
For simple $\lambda_0^{(l)}$, the corresponding eigenvector $e^{(l)}$ is even or odd; for double $\lambda_0^{(l)}=\lambda_0^{(l+1)}$, we can choose $e^{(l)}$ even and $e^{(l+1)}$ odd; then $\{e^{(l)}\}_{l=1}^{n-1}$
forms the basis in $\mathbb{R}^{n-1}$.
\par 2. Similar properties are valid for the eigenpairs of the first problem \eqref{eq:eig_elem} with the exception of one simple zero eigenvalue.
\end{lemma}

One can check by the direct computation that all the eigenvalues in $S_n$ and $\tilde{S}_n$ are simple at least for $1\leqslant n\leqslant 9$, see \cite{ZZ12}.
For low $n$, one can find $S_n$ and $\tilde{S}_n$ exactly, in particular,
$\tilde{S}_2=\{2.5\}$,
$\tilde{S}_3=\{2.5,10.5\}$,
$\tilde{S}_4=\{14\pm\sqrt{133},10.5\}$ and
$\tilde{S}_5=\{14\pm\sqrt{133},30\pm9\sqrt{5}\}$.

\par We choose $\{e^{(l)}\}_{l=1}^{n-1}$ as in Lemma \ref{prop1} using scaling $\widetilde{C}e^{(l)}\cdot e^{(l)}=1$.
\begin{lemma}
\label{lem1}
Let $\widetilde{G}(\lambda)p=-g(\lambda)$, see \eqref{eq:matrAC}, where  $\lambda\not\in\tilde{S}_n$.
Then the following formulas hold
\begin{gather*}
 p=\sum\nolimits_{l=1}^{n-1}\frac{a^{(l)}-\lambda c^{(l)}}{\lambda-\lambda_0^{(l)}}e^{(l)}
  =\sum\nolimits_{l=1}^{n-1}\frac{a^{(l)}-\lambda_0^{(l)} c^{(l)}}{\lambda-\lambda_0^{(l)}}e^{(l)}-\widetilde{C}^{-1}c.
\label{eq:formp}
\end{gather*}
Here $\{a^{(l)}\}_{l=1}^{n-1}$ and $\{c^{(l)}\}_{l=1}^{n-1}$ are the expansion coefficients of the vectors $a$ and $c$, see \eqref{eq:matrAC}, with
respect to the basis $\{\widetilde{C}e^{(l)}\}_{l=1}^{n-1}$, for example,
$c=\sum\nolimits_{l=1}^{n-1}c^{(l)}\widetilde{C}e^{(l)}$ with $c^{(l)}=c\cdot e^{(l)}$.
\end{lemma}

\smallskip\par 2. Below we need to assume that all the eigenvalues in both $S_n$ and $\tilde{S}_n$ are simple for considered $n$.
We introduce the auxiliary equation
\[
 \widehat{\gamma}(\lambda)\equiv -(g_0-g\cdot\widetilde{G}^{-1}g)(\lambda)/(g_n-\check{g}\cdot\widetilde{G}^{-1}g)(\lambda)=\theta
\]
with the parameter $\theta$, see \cite{ZZ12}.
Owing to Lemma \ref{lem1} this equation can be rewritten as
\begin{gather}
 a_0-\lambda c_0+\sum\nolimits_{l=1}^{n-1}\frac{(a^{(l)}-\lambda c^{(l)})^2}{\lambda-\lambda_0^{(l)}}
 =-\theta\Big(a_n-\lambda c_n+\sum\nolimits_{l=1}^{n-1}\frac{(\check{a}^{(l)}-\lambda \check{c}^{(l)})(a^{(l)}-\lambda c^{(l)})}{\lambda-\lambda_0^{(l)}}\Big).
\label{eq:theta}
\end{gather}
Its solving is equivalent to finding the roots of a polynomial having at most $n$th degree.
Here
$\check{a}^{(l)}=\check{a}\cdot e^{(l)}$ and $\check{c}^{(l)}=\check{c}\cdot e^{(l)}$.
Moreover, for $2\leqslant n\leqslant 9$ computations help to confirm that the vectors $e^{(l)}$ are even and odd respectively for odd and even $l$; therefore $\check{a}^{(l)}=(-1)^la^{(l)}$ and $\check{c}^{(l)}=(-1)^lc^{(l)}$, $l=\overline{1,n-1}$.
\par We define the simplest inner product in $S_K^{(n)}$ and the squared $\mathcal{C}$-norm
\[
 (y,v)_{S_K^{(n)}}:=\sum\nolimits_{j=1}^{K-1}y_jv_j+\sum\nolimits_{j=1}^{K}y_{j-1/2}\cdot v_{j-1/2},\ \
 \|v\|_{\mathcal{C}}^2:=(\mathcal{C}v,v)_{S_K^{(n)}}.
\]
Next theorem presents eigenvalues and eigenvectors of problem \eqref{eq:eig_glob}.
\begin{theorem}
\label{th:eigpares}
1. The spectrum of problem \eqref{eq:eig_glob} consists in $\tilde{S}_n$
and the numbers $\bigl\{\lambda_k^{(l)}\bigr\}_{l=1}^n\not\subset \tilde{S}_n$ that are all $n$ (and all positive real) solutions to equation \eqref{eq:theta} with $\theta=\theta_k:=\cos\frac{\pi k}{K}$ for $k=\overline{1,K-1}$
and are different for fixed $k$.
\par 2. To the eigenvalue $\lambda_0^{(l)}$, the following eigenvector corresponds
\begin{gather*}
 s_{0,j}^{(l)}=0,\,\ j=\overline{1,K-1},\ \
 s_{0,j-1/2}^{(l)}=(-P)^{j-1}e^{(l)},\,\ j=\overline{1,K},
\label{eq:eig vec1}
\end{gather*}
for $l=\overline{1,n-1}$. Here $(-P)^{j-1}e=(-1)^{j-1}e$ for even $e$, $(-P)^{j-1}e=e$ for odd $e$.
\par 3. To the eigenvalue $\lambda_k^{(l)}$, the following eigenvector corresponds
\begin{gather*}
s_{k,j}^{(l)}=\sin\frac{\pi kj}{K},\,\ j=\overline{1,K-1},\ \
 s_{k,j-1/2}^{(l)}=p_k^{(l)}\sin\frac{\pi k(j-1)}{K}+\check{p}_k^{(l)}\sin\frac{\pi kj}{K},\,\ j=\overline{1,K},
\label{eq:eig vec2}
\end{gather*}
where $p_k^{(l)}\in\mathbb{R}^{n-1}$ is the solution to non-degenerate algebraic system
$\widetilde{G}\bigl(\lambda^{(l)}_k\bigr)p_k^{(l)}=-g\bigl(\lambda^{(l)}_k\bigr)$, for $k=\overline{1,K-1}$, $l=\overline{1,n}$.

\par 4. The introduced eigenvectors are $C$-orthogonal, i.e. $(Cs_{k}^{(l)},s_{\tilde{k}}^{(\tilde{l})})_{S_K^{(n)}}=0$ for any $k,\tilde{k}\in\overline{0,K-1}$, $l\in\overline{1,n-\delta_{k0}}$ and $\tilde{l}\in\overline{1,n-\delta_{\tilde{k}0}}$ such that $k\neq\tilde{k}$ and/or $l\neq\tilde{l}$.
\par They form the basis in $S_K^{(n)}$, i.e. any $w\in S_K^{(n)}$ can be uniquely expanded as
\begin{gather}
 w=\sum\nolimits_{l=1}^{n-1}w_{0l}s_0^{(l)}+\sum\nolimits_{k=1}^{K-1}\sum\nolimits_{l=1}^n w_{kl}s_k^{(l)}.
\label{eq:decomp}
\end{gather}
\end{theorem}

\par Notice that:
(1) the vectors $s_0^{(l)}$ are used only to describe the algorithm, and only the vectors $e^{(l)}$ are applied in its implementation;
(2) $s_{k,j}^{(l)}$ are independent on $l$;
(3) the vectors $p_k^{(l)}$ can also be computed owing to Lemma \ref{lem1}.

\smallskip\par 3. We call the calculation of $w\in S_K^{(n)}$ by the coefficients $w_{kl}$ of the expansion \eqref{eq:decomp} as \textit{the inverse $F_n$-transform} and the calculation of the coefficients $w_{kl}$ by $w\in S_K^{(n)}$ as \textit{the direct $F_n$-transform}.
Let us describe their fast FFT-based implementation.
\begin{theorem}
\label{th:eigpares}
1. The inverse $F_n$-transform can be implemented according to the following formulas
\begin{gather*}
 w_j=\sum\nolimits_{k=1}^{K-1}\Big(\sum\nolimits_{l=1}^n w_{kl}\Big)\sin\frac{\pi kj}{K},\ \ j=\overline{1,K-1},
\label{eq:decomp_1}\\
 w_{j-1/2}=(-P)^{j-1}\sum\nolimits_{l=1}^{n-1}w_{0l}e^{(l)}
\\
 +2\sum\nolimits_{k=1}^{K-1}d_{k,e}\cos\frac{\pi k}{2K}\sin\frac{\pi k(j-1/2)}{K}
 -2\sum\nolimits_{k=1}^{K-1}d_{k,o}\sin\frac{\pi k}{2K}\cos\frac{\pi k(j-1/2)}{K},\ \ j=\overline{1,K},
\label{eq:decomp_2}
\end{gather*}
where $d_{k,e}$ and $d_{k,o}$ are respectively even and odd components of the vectors
$d_k:=\sum\nolimits_{l=1}^n w_{kl}p_k^{(l)}$.
Note that $(-P)^{j-1}e=e$ for odd $j$ and $(-P)^{j-1}e=-\check{e}$ for even $j$ for any $e\in\mathbb{R}^{n-1}$.
\par The collection $\{w_j\}_{j=1}^{K-1}$ can be computed by the standard inverse FFT with respect to sines.
The collection  $\{w_{j-1/2}\}_{j=1}^K$ can be computed by $n-1$ modified inverse FFT related to  the centers of elements in the amount of $[n/2]$ with respect to sines and $[(n-1)/2]$  with respect to cosines using extensions $d_{K,e}:=0$ and $d_{0,o}:=0$, see algorithms DST-I, DST-III and DCT-III in \cite{BRY07}.
\par 2. The direct $F_n$-transform can be implemented starting from the standard formulas
\begin{gather*}
 w_{kl}=(\mathcal{C}w,s_{k}^{(l)})_{S_K^{(n)}}/\|s_k^{(l)}\|_{\mathcal{C}}^2.
\label{eq:coef decomp_2}
\end{gather*}
Here, first, for $k=0$, $l=\overline{1,n-1}$, we have
\begin{gather*}
  (\mathcal{C}w,s_0^{(l)})_{S_K^{(n)}}
 =\Big(\widetilde{C}\sum\nolimits_{j=1}^K(-P)^{j-1}w_{j-1/2}\Big)\cdot e^{(l)}, \ \
 \|s_0^{(l)}\|_{\mathcal{C}}^2=K.
\label{eq:norms eigf}
\end{gather*}
\par Second, for $k=\overline{1,K-1}$, $l=\overline{1,n}$ and $y:=Cw$,
we have
\begin{gather*}
 (y,s_k^{(l)})_{S_K^{(n)}}=\sum\nolimits_{j=1}^{K-1}y_j\sin\frac{\pi kj}{K}
\nonumber\\
 +p_{k,e}^{(l)}\cdot\sum\nolimits_{j=1}^{K-1}(y_{j-1/2}+y_{j+1/2})_e\sin\frac{\pi kj}{K}
 +p_{k,o}^{(l)}\cdot\sum\nolimits_{j=1}^{K-1}(y_{j+1/2}-y_{j-1/2})_o\sin\frac{\pi kj}{K},
\label{eq:inv_fn}\\
 \|s_k^{(l)}\|_{\mathcal{C}}^2=K(b_{kl,0}+b_{kl,n}\theta_k),\,\
 b_{kl,0}=c_0+(\widetilde{C}p_k^{(l)}+2c)\cdot p_k^{(l)},\,\
 b_{kl,n}=c_n+(\widetilde{C}p_k^{(l)}+2c)\cdot \check{p}_k^{(l)}.
\nonumber
\end{gather*}
The collection of all these coefficients can be computed using $n$ standard direct FFTs  with respect to sines.
\end{theorem}

\smallskip\par 4. Now we consider in detail solving of the $N$-dimensional boundary value problem
\begin{gather}
 -\Delta u+\alpha u=f\ \ \text{в}\ \ \Omega=(0,X_1)\times\ldots\times(0,X_N),\ \ u|_{\partial\Omega}=0,
\label{eq:diff bvp pr2}
\end{gather}
where $\Delta$ is the Laplace operator and $\alpha=\textrm{const}$; for simplicity, let $\alpha>-\pi^2\bigl(X_1^{-2}+\ldots+X_N^{-2}\bigr)$.
\par We introduce the space $H_{h_1}^{(n_1)}[0,X_1]\otimes\ldots\otimes H_{h_N}^{(n_N)}[0,X_N]$ of the piecewise-poly\-no\-mial in
$\overline{\Omega}$ functions, where  $h_i=X_i/K_i$ and $n_i\geq 2$, $i=\overline{1,N}$.
Let $\mathbf{K}=(K_1,\ldots,K_N)$ and $\mathbf{n}=(n_1,\ldots,n_N)$.
\par We define the space $S_{\mathbf{K}}^{(\mathbf{n})}=S_{K_1}^{(n_1)}\otimes\ldots\otimes S_{K_N}^{(n_N)}$ of vector functions.
Similarly to the 1D case, there is the natural isomorphism between
functions in $H_{h_1}^{(n_1)}[0,X_1]\otimes\ldots\otimes H_{h_N}^{(n_N)}[0,X_N]$ and vectors in $S_{\mathbf{K}}^{(\mathbf{n})}$.
\par The FEM dicretization of problem \eqref{eq:diff bvp pr2} can be written in the following operator form
\begin{gather}
 (4h_1^{-2}\mathcal{A}_1\mathcal{C}_2\ldots\mathcal{C}_N+\ldots+4h_N^{-2}\mathcal{A}_N\mathcal{C}_1\ldots\mathcal{C}_{N-1})v
 +\alpha \mathcal{C}_1\ldots\mathcal{C}_Nv=f^h,\ \ v\in S_{\mathbf{K}}^{(\mathbf{n})},
\label{eq:bvp glob 2}
\end{gather}
where $\mathcal{A}_i$ and $\mathcal{C}_i$ are versions of the above defined operators $\mathcal{A}$ and $\mathcal{C}$ acting in variable  $x_i$ (depending on $K_i$ and $n_i$), $i=\overline{1,N}$, and $f^h\in S_{\mathbf{K}}^{(\mathbf{n})}$ is the FEM average of $f$.
Remind that the general case $u|_{\partial\Omega}=b$ in \eqref{eq:diff bvp pr2} could be covered by reducing to \eqref{eq:bvp glob 2} with the modified $f^h$ depending on $b^h$ (the FEM average of $b$).
\par To compute its solution, the $F_n$-transforms from Theorem \ref{th:eigpares} can be applied twofold.
\smallskip\par (a) Let the vector $\varphi^h\in S_{\mathbf{K}}^{(\mathbf{n})}$ be the solution to the auxiliary algebraic problem
$\mathcal{C}_1\ldots\mathcal{C}_N\varphi^h=f^h$
with the splitting operator (the product of operators acting in $x_1,\ldots,x_N$),
i.e. formally $\varphi^h=\mathcal{C}_1^{-1}\ldots\mathcal{C}_N^{-1}f^h$.
We consider the multiple expansion of $\varphi^h\in S_{\mathbf{K}}^{(\mathbf{n})}$ like \eqref{eq:decomp}
\begin{gather}
 \varphi^h
 =\sum\nolimits_{i=1}^N\sum\nolimits_{k_i=0}^{K_i-1}\,\sum\nolimits_{l_i=1}^{n_i-\delta_{k_i0}}
 \varphi^h_{k_1l_1,\ldots,k_Nl_N}s_{1,\,k_1}^{(l_1)}\ldots s_{N,\,k_N}^{(l_N)}.
\label{eq:decomp phi 2}
\end{gather}
Then the expansion of the solution has the following form
\begin{gather}
 v=\sum\nolimits_{i=1}^N\sum\nolimits_{k_i=0}^{K_i-1}\,\sum\nolimits_{l_i=1}^{n_i-\delta_{k_i0}}
 \frac{\varphi^h_{k_1l_1,\ldots,k_Nl_N}}{4h_1^{-2}\lambda_{1,\,k_1}^{(l_1)}+\ldots+4h_N^{-2}\lambda_{m,\,k_N}^{(l_N)}+\alpha}
 s_{1,\,k_1}^{(l_1)}\ldots s_{N,\,k_N}^{(l_N)}.
\label{eq:decomp v 2}
\end{gather}
Here $\bigl\{\lambda_{i,k_i}^{(l_i)},s_{i,k_i}^{(l_i)}\bigr\}$ are versions of the above defined eigenpairs
$\bigl\{\lambda_{k}^{(l)},s_{k}^{(l)}\bigr\}$ with respect to $x_i$.
\par The steps of the algorithm (a) are rather standard:
\par (1) solving the auxiliary problem $\mathcal{C}_1\ldots\mathcal{C}_N\varphi^h=f^h$
for $\varphi^h$ (that is reduced to the sequential solving of the 1D problems in $x_1$ with the matrix $\mathcal{C}_1$,..., $x_N$ with the matrix  $\mathcal{C}_N$);
\par (2) finding the coefficients of expansion \eqref{eq:decomp phi 2} for $\varphi^h$ (by the direct $F_n$-transforms in $x_1$,..., $x_N$);
\par (3) finding $v$ by the coefficients of its expansion \eqref{eq:decomp v 2} (by the inverse $F_n$-transforms in $x_1$,..., $x_N$).

\par (b) Let the vector $\varphi^h\in S_{\mathbf{K}}^{(\mathbf{n})}$ be the solution to the auxiliary $(m-1)$D problem
$\mathcal{C}_2\ldots\mathcal{C}_N\varphi^h=f^h$
in $x_2$,..., $x_N$, i.e. formally $\varphi^h=\mathcal{C}_2^{-1}\ldots\mathcal{C}_N^{-1} f^h$.
We consider the expansion of $\varphi^h$ like \eqref{eq:decomp} in $x_2$,..., $x_N$, i.e.
\begin{gather}
 \varphi^h=
 \sum\nolimits_{i=2}^N\sum\nolimits_{k_i=0}^{K_i-1}\,
 \sum\nolimits_{l_i=1}^{n_i-\delta_{k_i0}}\varphi^h_{k_2l_2,\ldots,k_Nl_N}s_{2,\,k_2}^{(l_2)}\ldots s_{N,\,k_N}^{(l_N)},
\label{eq:decomp phi}
\end{gather}
now with the coefficients $\varphi^h_{k_2l_2,\ldots,k_Nl_N}\in S_{K_1}^{(n_1)}$.
Then the coefficients $v_{kl}\in S_{K_1}^{(n_1)}$ in the similar expansion of the solution $v\in S_{\mathbf{K}}^{(\mathbf{n})}$
\begin{gather}
 v=
 \sum\nolimits_{i=2}^N\sum\nolimits_{k_i=0}^{K_i-1}\,\sum\nolimits_{l_i=1}^{n_i-\delta_{k_i0}}
 v_{k_2l_2,\ldots,k_Nl_N}s_{2,\,k_2}^{(l_2)}\ldots s_{N,\,k_N}^{(l_N)},
\label{eq:decomp v}
\end{gather}
serve as the solutions to 1D problems in $x_1$
\begin{gather}
 \bigl[4h_1^{-2}\mathcal{A}_1
 +(4h_2^{-2}\lambda_{k_2}^{(l_2)}+\ldots+4h_N^{-2}\lambda_{k_N}^{(l_N)}+\alpha)\mathcal{C}_1\bigr]v_{k_2l_2,\ldots,k_Nl_N}
 =\varphi^h_{k_2l_2,\ldots,k_Nl_N}.
\label{eq:sys x1 2}
\end{gather}
Their matrices are symmetric and positive definite.
Of course, the simpler case $n_1=1$ is acceptable too.
\par The steps of the algorithm (b) are rather standard as well:
\par (1) solving the auxiliary problem $\mathcal{C}_2\ldots\mathcal{C}_N\varphi^h=f^h$
for $\varphi^h$;
\par (2)  finding the coefficients of the expansion \eqref{eq:decomp phi} for $\varphi^h$ (by the direct $F_n$-transforms in $x_2$,..., $x_N$);
\par (3) solving the collection of the 1D problems \eqref{eq:sys x1 2} for the coefficients of the expansion of $v$;
\par (4) finding $v$ by the coefficients of its expansion  \eqref{eq:decomp v} (by the inverse $F_n$-transforms in $x_2$,..., $x_N$).

\par Implementing algorithms (a) and (b) costs respectively $O(K\log_2K)$ and
$O\bigl(K\log_2(K_2\ldots K_N)\bigr)$ arithmetic operations with $K=K_1\ldots K_N$.

They can be applied to solve various time-dependent PDEs such as the heat, wave or Schr\"{o}dinger's equations since usually their implicit time discretizations lead to problems like \eqref{eq:bvp glob 2} at the upper time level.

\par Moreover, algorithm (b) is directly extended to the case of more general equations than in \eqref{eq:diff bvp pr2} with the coefficients depending on $x_1$, various boundary conditions for $x_1=0,X_1$ and the nonuniform mesh in $x_1$ \cite{SN78}.
It can also be applied to reduce 3D problems in a cylindrical domain to a collection of independent 2D problems in the cylinder base.

\smallskip\par 5. Both algorithms (a) and (b) are well-behaved in the numerical experiments.
We choose problem \eqref{eq:diff bvp pr2} for $N=2$, $\alpha=1$ and $X_1=X_2=1$, with the exact solution
$u(x_1,x_2):=\sin(\pi x_1)\sin(\pi x_2)(x_1+x_2-1)$ and take $K_1=K_2=K$ and $n_1=n_2=n$.
The errors for algorithm (a) in the uniform norm are given in Table \ref{tab:EX42:2D:ErrorCabs} in dependence on $K=8,16,...,1024$, for $n=\overline{1,9}$.
We emphasize that there is almost no impact of the round-off errors as $K$ and $n$ grows.
Here the multiple Gauss quadrature formulas with $n+1$ nodes in $x_1$ and $x_2$ were applied to compute $f^h$, and the eigenvalues of the 1D problems were computed with the quadruple precision (using Mathematica) to improve the stability with respect to round-off errors.
\par
\begin{figure}[htbp]\centering{
    \begin{minipage}[h]{0.49\linewidth}
        \center{\includegraphics[width=1\linewidth]{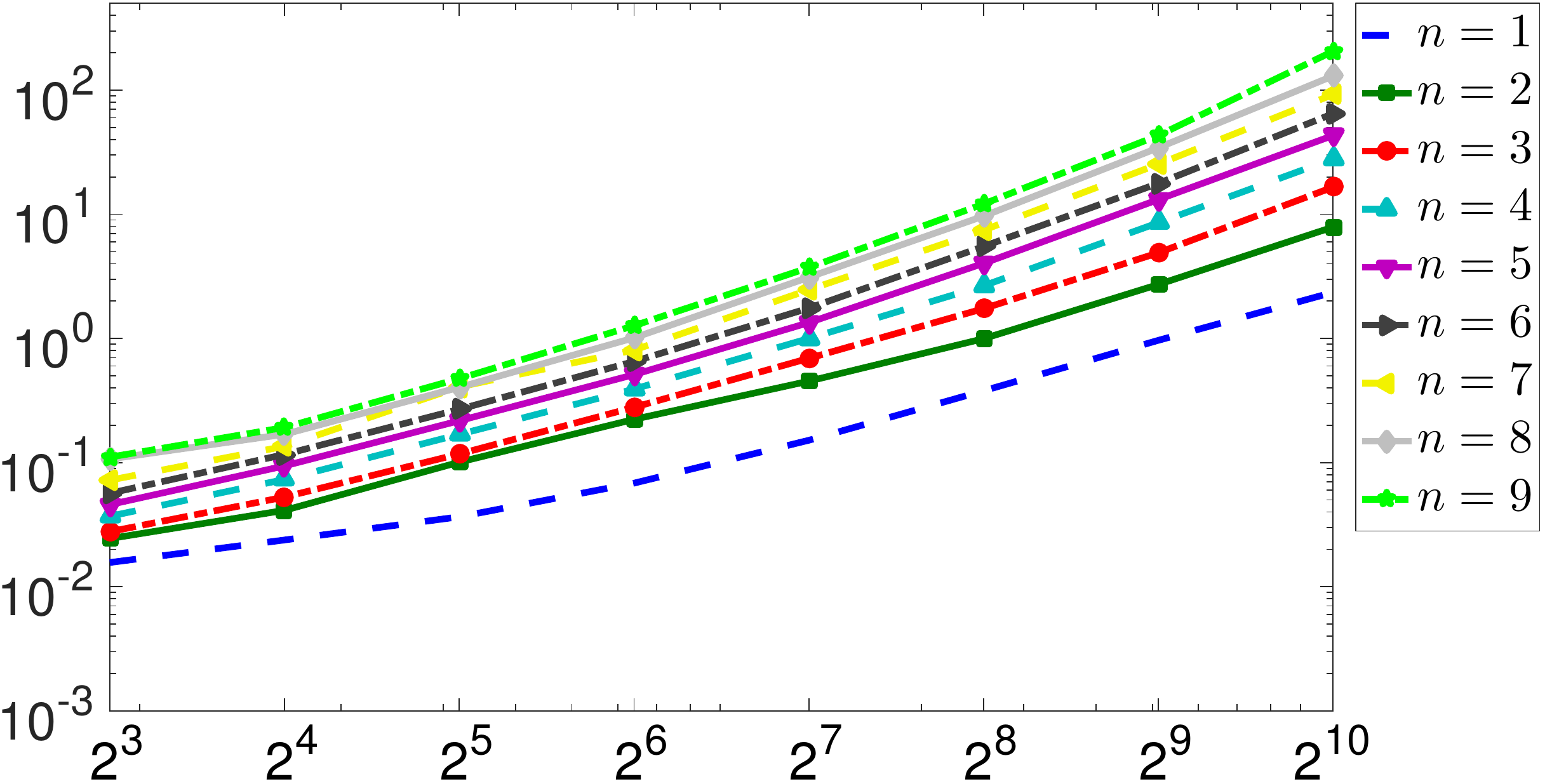}}
    \end{minipage}
    \begin{minipage}[h]{0.49\linewidth}
        \center{\includegraphics[width=1\linewidth]{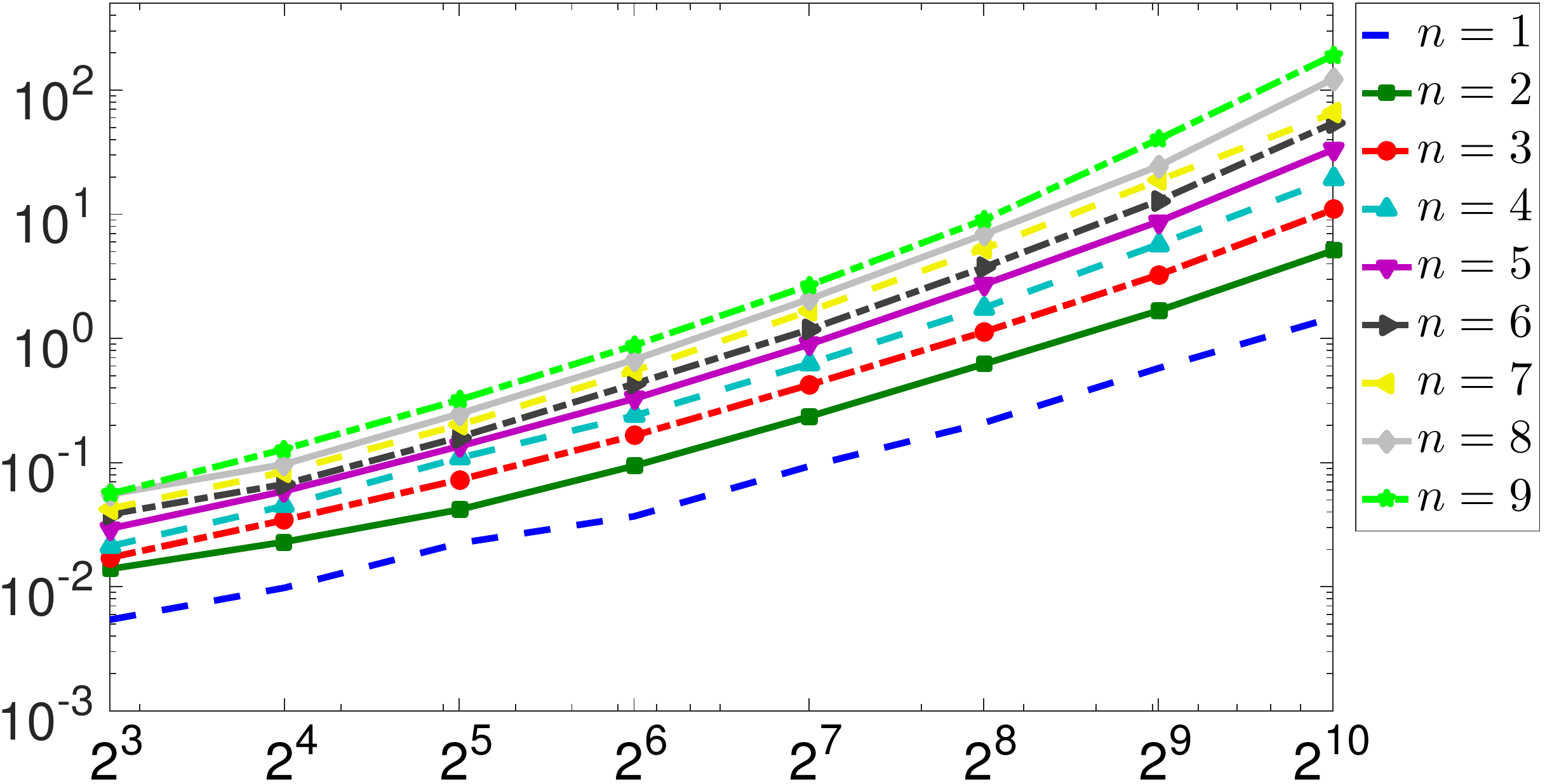}}
    \end{minipage}
    }\caption{\small{The execution time (in seconds) for algorithms (a) (left) and (b) (right)}}
\label{fig:TIME:2D:ASUS}
\end{figure}

 In Fig. \ref{fig:TIME:2D:ASUS} we present the execution time for the same $K$ and $n$, using our codes in Matlab R2016a for both algorithms.
The ordinary laptop with Intel Core i3-2350M CPU 2.3 GHz, 4 Gb, Win 7 x64 on board was applied.
Including the case $n=1$ allows us to compare the original well-known algorithms with the above suggested new algorithms for higher $n$.
Notice the rather close to linear behavior of time in $K$ and its mild monotone growth in $n$.
Specify that system \eqref{eq:bvp glob 2} contains $(Kn-1)^2$ unknowns.
For $K=1024$ and $n=9$, this is almost $85\cdot10^6$ unknowns but only less than 2 min is required for solving.

\smallskip\par\textbf{Acknowledgement.}
The study has been funded
by the RFBR, grant № 16-01-00048.
\begin{table}\centering{
    \begin {tabular}{r<{\pgfplotstableresetcolortbloverhangright }@{}l<{\pgfplotstableresetcolortbloverhangleft }r<{\pgfplotstableresetcolortbloverhangright }@{}l<{\pgfplotstableresetcolortbloverhangleft }r<{\pgfplotstableresetcolortbloverhangright }@{}l<{\pgfplotstableresetcolortbloverhangleft }r<{\pgfplotstableresetcolortbloverhangright }@{}l<{\pgfplotstableresetcolortbloverhangleft }r<{\pgfplotstableresetcolortbloverhangright }@{}l<{\pgfplotstableresetcolortbloverhangleft }r<{\pgfplotstableresetcolortbloverhangright }@{}l<{\pgfplotstableresetcolortbloverhangleft }r<{\pgfplotstableresetcolortbloverhangright }@{}l<{\pgfplotstableresetcolortbloverhangleft }r<{\pgfplotstableresetcolortbloverhangright }@{}l<{\pgfplotstableresetcolortbloverhangleft }r<{\pgfplotstableresetcolortbloverhangright }@{}l<{\pgfplotstableresetcolortbloverhangleft }r<{\pgfplotstableresetcolortbloverhangright }@{}l<{\pgfplotstableresetcolortbloverhangleft }}%
\toprule \multicolumn {2}{c}{\phantom{X}$\kappa$\phantom{X}}&\multicolumn {2}{c}{\phantom{x}$n=1$\phantom{x}}&\multicolumn {2}{c}{\phantom{x}$n=2$\phantom{x}}&\multicolumn {2}{c}{\phantom{x}$n=3$\phantom{x}}&\multicolumn {2}{c}{\phantom{x}$n=4$\phantom{x}}&\multicolumn {2}{c}{\phantom{x}$n=5$\phantom{x}}&\multicolumn {2}{c}{\phantom{x}$n=6$\phantom{x}}&\multicolumn {2}{c}{\phantom{x}$n=7$\phantom{x}}&\multicolumn {2}{c}{\phantom{x}$n=8$\phantom{x}}&\multicolumn {2}{c}{\phantom{x}$n=9$\phantom{x}}\\\midrule %
%$2$&$$&$$&$$-&$3$&$.7E^{-3}$&$1$&$.8E^{-3}$&$2$&$.9E^{-4}$&$3$&$.7E^{-5}$&$2$&$.0E^{-6}$&$1$&$.4E^{-7}$&$5$&$.4E^{-9}$&$2$&$.9E^{-10}$\\%
%$4$&$$&$2$&$.3E^{-2}$&$3$&$.3E^{-4}$&$1$&$.3E^{-4}$&$1$&$.4E^{-5}$&$5$&$.1E^{-7}$&$2$&$.3E^{-8}$&$4$&$.8E^{-10}$&$1$&$.6E^{-11}$&$2$&$.6E^{-13}$\\%
$3$&$$&$5$&$.4E^{-3}$&$1$&$.9E^{-5}$&$8$&$.6E^{-6}$&$5$&$.0E^{-7}$&$8$&$.2E^{-9}$&$2$&$.0E^{-10}$&$1$&$.9E^{-12}$&$3$&$.5E^{-14}$&$1$&$.2E^{-14}$\\%
$4$&$$&$1$&$.4E^{-3}$&$1$&$.2E^{-6}$&$5$&$.4E^{-7}$&$1$&$.6E^{-8}$&$1$&$.3E^{-10}$&$1$&$.6E^{-12}$&$8$&$.5E^{-15}$&$1$&$.0E^{-14}$&$4$&$.5E^{-14}$\\%
$5$&$$&$3$&$.5E^{-4}$&$7$&$.7E^{-8}$&$3$&$.4E^{-8}$&$5$&$.1E^{-10}$&$1$&$.9E^{-12}$&$1$&$.3E^{-14}$&$9$&$.3E^{-15}$&$1$&$.3E^{-14}$&$6$&$.4E^{-14}$\\%
$6$&$$&$8$&$.7E^{-5}$&$4$&$.8E^{-9}$&$2$&$.1E^{-9}$&$1$&$.6E^{-11}$&$3$&$.0E^{-14}$&$6$&$.1E^{-16}$&$5$&$.6E^{-15}$&$5$&$.2E^{-15}$&$1$&$.7E^{-14}$\\%
$7$&$$&$2$&$.2E^{-5}$&$3$&$.0E^{-10}$&$1$&$.3E^{-10}$&$4$&$.9E^{-13}$&$2$&$.4E^{-15}$&$9$&$.4E^{-16}$&$6$&$.4E^{-15}$&$7$&$.4E^{-15}$&$4$&$.5E^{-15}$\\%
$8$&$$&$5$&$.4E^{-6}$&$1$&$.9E^{-11}$&$8$&$.3E^{-12}$&$1$&$.6E^{-14}$&$1$&$.9E^{-15}$&$1$&$.6E^{-15}$&$7$&$.3E^{-15}$&$9$&$.2E^{-15}$&$1$&$.6E^{-14}$\\%
$9$&$$&$1$&$.4E^{-6}$&$1$&$.2E^{-12}$&$5$&$.2E^{-13}$&$6$&$.9E^{-16}$&$1$&$.3E^{-15}$&$1$&$.4E^{-15}$&$2$&$.7E^{-15}$&$8$&$.5E^{-15}$&$2$&$.8E^{-14}$\\%
$10$&$$&$3$&$.4E^{-7}$&$7$&$.4E^{-14}$&$3$&$.3E^{-14}$&$6$&$.1E^{-16}$&$6$&$.9E^{-16}$&$1$&$.7E^{-15}$&$2$&$.2E^{-15}$&$1$&$.6E^{-14}$&$3$&$.4E^{-14}$\\\bottomrule %
\end {tabular}%
}
\caption{Errors in the uniform norm in dependence on $K=2^\kappa=8,16,...,1024$ and $n=\overline{1,9}$
\label{tab:EX42:2D:ErrorCabs}}
\end{table}

\end{document}